\documentstyle[12pt, amssymbols]{article}

\def\nat{\cal N}
\newtheorem{theoremfoo}{Theorem}[section] 
\newenvironment{theorem}{\pagebreak[1]\begin{theoremfoo}}{\end{theoremfoo}}
\newtheorem{propositionfoo}[theoremfoo]{Proposition}

\newtheorem{lemmafoo}[theoremfoo]{Lemma}

\newtheorem{corollaryfoo}[theoremfoo]{Corollary}
\newenvironment{corollary}{\pagebreak[1]\begin{corollaryfoo}}{\end{corollaryfoo}}
\newtheorem{dfntn}[theoremfoo]{Definition}
\newenvironment{definition}{\pagebreak[1]\begin{dfntn}\rm}{\end{dfntn}}
\newenvironment{proof}
    {\pagebreak[1]{\narrower\noindent {\bf Proof:\quad\nopagebreak}}}{\QED}
\newenvironment{sketch}
    {\pagebreak[1]{\narrower\noindent {\bf Proof sketch:\quad\nopagebreak}}}{\QED}

\newcommand{\yyskip}{\penalty-50\vskip 5pt plus 3pt minus 2pt}
\newcommand{\blackslug}{\hbox{\hskip 1pt
        \vrule width 4pt height 8pt depth 1.5pt\hskip 1pt}}
\newcommand{\QED}{{\penalty10000\parindent 0pt\penalty10000
        \hskip 8 pt\nolinebreak\blackslug\hfill\lower 8.5pt\null}
        \par\yyskip\pagebreak[1]}
\newtheorem{factfoo}[theoremfoo]{Fact}

\begin{document}
\title{Kolmogorov Complexity and Instance Complexity of Recursively Enumerable Sets} 
\date{}
\author{
{Martin Kummer\/}\thanks{
Institut f\"ur Logik,
Komplexit\"at und Deduktionssysteme,
Universit\"at Karlsruhe, D-76128 Karlsruhe, Germany,
(Email: {\tt kummer@ira.uka.de}).}
\\ {\small Universit\"at Karlsruhe}
}
\maketitle

\begin{abstract} We study in which way Kolmogorov complexity
and instance complexity affect properties of r.e.\ sets. 
We show that the well-known $2\log n$ upper bound on the Kolmogorov 
complexity of initial segments of r.e.\ sets is optimal and characterize
the T-degrees of r.e.\ sets which attain this bound.
The main part of the paper is concerned with instance
complexity of r.e.\ sets. We construct a nonrecursive r.e.\ set
with instance complexity logarithmic in the Kolmogorov complexity.
This refutes a conjecture of Ko, Orponen, Sch{\"o}ning, and Watanabe.
In the other extreme, we show that all wtt-complete set and
all Q-complete sets have infinitely many hard instances.
\end{abstract}



\section{Introduction}

Intuitively, Kolmogorov complexity measures the ``descriptional complexity'' of a string $x$. It is defined as 
the length of the shortest program that computes $x$ from the empty input. 
Accordingly, the Kolmogorov complexity of initial segments of a set $A$ 
is considered as a measure of the ``randomness'' of $A$. It is well-known that for r.e.\ sets 
the Kolmogorov complexity of initial segments of length $n$ is bounded by $2 \log n$. 
We show that this bound is optimal and characterize the Turing degrees of r.e.\ sets which attain
this bound as the array nonrecursive degrees of Downey, Jockusch, and Stob \cite{DJS90}.

Ko, Orponen, Sch{\"o}ning, and Watanabe \cite{K+86, O+94} have recently 
introduced the notion of {\it instance
complexity\/} as a measure of the complexity of individual instances of $A$. Informally,
$ic(x:A)$, the instance complexity of $x$ with respect to $A$, is the length of the
shortest total program which correctly computes $\chi_A(x)$ and does not make any
mistakes on other inputs, but it is permitted to output ``don't know'' answers.
It is easy to see that the Kolmogorov complexity of $x$ is an upper bound for the
instance complexity of $x$ (up to a constant). A set $A$ has {\it hard instances\/}
if for infinitely many $x$ the instance complexity of $x$ w.r.t.\ $A$ is at least as high as the
Kolmogorov complexity of $x$ (up to a constant which may depend on $A$), i.e., the trivial
upper bound is already optimal.

Orponen et al.\ conjectured in \cite{Or90, O+94} that every nonrecursive r.e.\ set has hard 
instances (``Instance Complexity Conjecture (ICC)''). Buhrmann and Orponen \cite{BO93} proved ICC for 
m-complete sets. Tromp \cite{Tr93} proved that the instance complexity of $x$ w.r.t.\  any nonrecursive set $A$
is infinitely often at least logarithmic in the Kolmogorov complexity of $x$. We construct
an r.e. nonrecursive set which attains this lower bound for all $x$. In particular,
this is a counterexample to ICC. On the positive side, we show that ICC holds for
wtt-complete sets, Q-complete sets, and hyperhypersimple sets. But ICC fails for a
T-complete set, since it fails for an effectively simple set. However,
ICC holds for all strongly effectively simple sets. We also investigate a
weak version of instance complexity, where programs may not halt instead of
giving ``don't know'' answers.

The  resource-bounded version of instance complexity is also well-studied;
we refer the reader to \cite{BO93, FK94, Ko92, O+94}.

\medskip
\noindent
{\it Notation and Definitions:\/}\\
The notation generally follows $\cite{LV93}$. For further recursion theoretic background
we refer the reader to \cite{Od89, So87}. For $p \in \{0,1\}^*$, $l(p)$ denotes the length
of $p$; $\lambda$ is the empty string. We use the special symbol $\perp$ to denote the
``don't know'' output. $\chi_A$ is the characteristic function of $A$. 
We identify $\nat$ and $\{0,1\}^*$ via the canonical correspondence as in \cite[p.\ 11]{LV93}.

\begin{definition} (Chaitin, Kolmogorov, Solomonoff)\\
For any partial recursive mapping 
$U : \{0,1\}^* \times \{0,1\}^* \rightarrow \{0,1\}^* \cup \{\perp\}$
and any $x \in \{0,1\}^*$ we define $C_U(x) = \min\{ l(p) : U(p,\lambda) = x\}$, the 
{\it Kolmogorov complexity\/} of $x$ in $U$.  If no such $p$ exists then $C_U(x) = \infty$.
\end{definition}

It is helpful to think of $U$ as an {\it interpreter\/} 
which takes a {\it program\/} $p$ and an input $z$
and produces the output $U(p,z)$

Instance complexity was introduced in \cite{K+86} 
in order to study the complexity
of single instances of a decision problem.

\begin{definition} (Ko, Orponen, Sch{\"o}ning, Watanabe, 1986)\\ 
Let $A \subseteq \{0,1\}^*$. A function $f : \{0,1\}^* \rightarrow \{0,1,\perp\}$ is called
$A$-consistent if $f(x) = \chi_A(x)\ \vee\ f(x) = \perp$, for all $x \in dom(f)$.
The {\it instance complexity\/} of $x$ with respect to $A$ in $U$ is defined as
\begin{tabbing}
\quad \quad $ic_U(x:A) = \min\{ l(p) :$ \= $\lambda z.\ U(p,z)$ is a {\it total\/} $A$-consistent function\\
\> such that $U(p,x) = \chi_A(x)\}$.
\end{tabbing}
If no such $p$ exists then $ic_U(x:A) = \infty$.
\end{definition}

If we drop in the definition of $ic_U$ the requirement that $\lambda z.\ U(p,z)$ is total,
then we obtain a weaker notion of instance complexity, which we denote by $\overline{ic}_U(x:A)$.
Note that $\overline{ic}_U(x:A) \leq ic_U(x:A)$ for all $x, A$.

It is well-known (see \cite{LV93}) that there exist ``optimal'' 
partial recursive functions $U$ such that,
for every partial recursive mapping $U'$, there is a constant $c$ with
$C_U(x) \leq C_{U'}(x) + c$, $ic_U(x:A) \leq ic_{U'}(x:A) + c$, and 
$\overline{ic}_U(x:A) \leq  \overline{ic}_{U'}(x:A) + c$, for all $x,A$.
 
For the following we fix an optimal mapping $U$ and write
$C(x)$, $ic(x:A)$, and $\overline{ic}(x:A)$, 
for $C_U(x)$,  $ic_U(x:A)$, and $\overline{ic}_U(x:A)$, respectively.
We also write $U_s(p,z)$ for the result, if any, after $s$ steps of computation
of $U$ with input $(p,z)$. $C^s(x)$ denotes the approximation to $C(x)$ after
$s$ steps of computation (i.e., with $U_s$ in place of $U$ in the definition of $C(x)$).
Clearly, $C^{s+1}(x) \leq C^s(x)$ and $C^t(x) = C(x)$ for all sufficiently large $t$.

The instance complexity of $x$ can be bounded by the Kolmogorov complexity of $x$
in the sense that for every set $A$ there is a constant $c$ such that
$ic(x:A) \leq C(x) + c$ for all $x$. Informally, $x$ is a hard instance of $A$ 
if this upper bound is also a lower bound. This was the motivation for the
following definition (which is independent of the choice of the optimal $U$).

\begin{definition} (Ko, Orponen, Sch{\"o}ning, Watanabe, 1986)\\
A set $A$ has {\it hard instances\/} if there is a constant $c$ such that 
$$ic(x:A) \geq C(x) -c \mbox{ \ for infinitely many } x.$$
If the condition holds with $\overline{ic}$ in place of $ic$, 
we say that $A$ has hard instances {\it with respect to $\overline{ic}$.}
\end{definition}

\noindent
{\it Remark:\/} 
The difference between $ic$ and $\overline{ic}$ is perhaps best explained by an example:

Suppose that $A$ is an r.e.\ set and we want to define a program $p$ such that it 
witnesses $ic(x:A) \leq |p|$ for all $x$ with $C(x) < n$. Since $p$ has to be total we 
have to define it for every input $z$ at some step $s$. If $z$ has already appeared in 
$A$ there is no problem, we set $U(p,z) = 1$. If $z$ has not yet appeared in $A$ and
$C^s(z) \geq  n$, we could try to define $U(p,z) = \perp$, but this can later become incorrect
if it turns out that $C(z) < n$. If we set $U(p,z) = 0$ and $z$ later appears in $A$,
then $p$ is also incorrect.

In the case of $\overline{ic}$ we have more freedom: We may leave  $U(p,z)$ 
undefined until $C^s(z) < n$ at some stage $s$. If this never happens, then
$U(p,z)$ is undefined and $C(z) \geq n$, which is fine. Still the second 
source of error remains: If $C^s(z) < n$ and $z$ has not yet appeared in 
$A$ at stage $s$, we have to define $U(p,z)$, and the best we can do is to set $U(p,z) = 0$. 
But this may later turn out to be incorrect.

\section{A version of Barzdin's Lemma}

In this section we consider the Kolmogorov complexity of initial segments of r.e.\ sets.
For $A \subseteq \nat$ and $n \in \nat$ we write $\chi_A \restriction n$ for the string
$\chi_A(0)\ldots\chi_A(n)$.

Let us first recall what was previously known.
The conditional complexity of a string $\sigma$ of length $n$ is
defined as $C(\sigma | n) = \min\{l(p) : U(p,n) = \sigma \}$.
We write $C(\chi_A | n)$ for $C(\chi_A \restriction (n-1) | n)$.
Barzdin (\cite{Ba68}, see \cite[Theorem 2.18]{LV93}) characterized the worst case
of the conditional complexity for initial segments of r.e.\ sets:

\begin{itemize}
\item For every r.e.\ set $A$ there is a constant $c$ such that for all $n$:\\
      $C(\chi_A | n) \leq \log n +c$.
\item There is an r.e.\ set $A$ such that $C(\chi_A | n) \geq \log n$ for all $n$.
\end{itemize}

Now we look at the standard Kolmogorov complexity $C(\chi_A \restriction n)$. 
Utilizing a result of Meyer \cite[p.\ 525]{Lo69},
Chaitin proved that, if there is constant $c$ such that  for all $n$,
$C(\chi_A \restriction n) \leq \log n + c$, then $A$ is recursive  
\cite[Theorem 6]{Ch76}, \cite[Exercise 2.43]{LV93}.

For every r.e.\ set $A$ there is a constant $c$ such that 
$C(\chi_A \restriction n) \leq 2 \log n + c $ for all $n$
(see \cite[Exercise 2.59]{LV93}).
On the other hand, there is no r.e.\ set $A$ such that $C(\chi_A \restriction n) \geq 2\log n - O(1)$
for almost all $n$. This follows from the argument in  \cite[Exercise 2.58]{LV93}.

In \cite[Exercise 2.59]{LV93} it is stated as an open question (attributed to Solovay)
whether the upper bound $2\log n$ is optimal. The following result shows that this is indeed the case. 
For ease of conversation, we say that $A$ is {\it complex\/} if there is a constant $c$
such that $C(\chi_A \restriction n)\ \geq\ 2 \log n - c$ for 
infinitely many $n \in \nat$.

\begin{theorem} \label{th:sol} There is an r.e.\ complex set.
\end{theorem}

\begin{proof} Let $t_0  = 0, t_{k+1} = 2^{t_k}$, and $I_k = (t_k,t_{k+1}]$, for all $k \geq 0$.
$(I_k)$ is a sequence of exponentially increasing half-open intervals.

Let $f(k) = \sum_{i = t_k+1}^{t_{k+1}} (i - t_k +1)$, $g(k) = \max\{l : 2^{l+1}-1 < f(k)\}$.
Note that $f(k) = \frac{1}{2} t^2_{k+1} - o(1)$ and $g(k) = 2 \log t_{k+1} -2 -o(1)$, for
$k \rightarrow \infty$.

We enumerate an r.e.\ set $A$ in steps as follows:
\begin{tabbing}
{\it Step $0$}: Let $A_0 = \emptyset$.\\
{\it Step $s+1$}: \= Let $A_{s+1} = A_s$. For $k = 0,\ldots,s$ do: 
If $C^s(\chi_{A_s} \restriction n) \leq g(k)$
for all $n \in I_k$\\
\> then enumerate $\min(\overline{A}_s \cap I_k)$ into $A_{s+1}$. $\Box$
\end{tabbing}
Let $A = \cup_{s \geq 0} A_s$.
Suppose for a contradiction that $C(\chi_{A} \restriction n) \leq g(k)$ for all $n \in I_k$.
Then we eventually enumerate every $n \in I_k$ into $A$. Note that for fixed $n$ there are at
least $n-t_k+1$ different strings $\sigma = \chi_{A_s} \restriction n$ with $l(\sigma) = n+1$
and $C(\sigma) \leq g(k)$. (The suffix of $\chi_{A_s} \restriction n$ runs through $1^x0^{n-t_k-x}$
for $x = 0,\ldots,n-t_k$.) Thus, there are at least $f(k)$ many different strings which all
have Kolmogorov complexity at most $g(k)$. This contradicts the definition of $g(k)$.

So for every $k$ there exists $n \in I_k$ with $C(\chi_A \restriction n) > g(k)$, i.e.,
$C(\chi_A \restriction n) > g(k) \geq 2 \log n - 2 - o(1)$. Thus, $A$ is complex.
\end{proof}

We now characterize the degrees of r.e.\ complex sets. Downey, Jockusch, and Stob \cite{DJS90}
introduced the notion of an {\it array nonrecursive\/} set. This captures precisely those
r.e.\ sets that arise in multiple permitting arguments. In \cite{DJS90} several other natural
characterizations of this degree class are given.

An r.e.\ set $A$ is called {\it array nonrecursive\/} with respect to $\{F_k\}_{k \in \nat}$ if
$$(\forall e)(\exists^{\infty} k) [ W_e \cap F_k = A \cap F_k ].$$
Here $\{F_k\}_{k \in \nat}$ denotes a very strong array. This means that $\{F_k\}_{k \in \nat}$ 
is a  strong array of pairwise disjoint sets which partition $\nat$ and satisfy$|F_k| < |F_{k+1}|$
for all $k \in \nat$.

An r.e.\ set is {\it array nonrecursive\/} if it is    
{\it array nonrecursive\/} with respect to some very strong array 
$\{F_k\}_{k \in \nat}$. A degree is called {\it array nonrecursive\/}
if it contains an r.e.\ array nonrecursive set. Not every r.e.\ nonrecursive
degree is array nonrecursive \cite[Theorem 2.10]{DJS90}.

\begin{theorem} \label{th:compdeg} The degrees containing an r.e.\ complex set 
coincide with the array nonrecursive degrees.
In addition, if $A$ is r.e.\ and not of array nonrecursive degree, 
then for every unbounded, nondecreasing, total recursive
function $f$ there is a constant $c$ such that 
$$C(\chi_A \restriction n) \leq \log n + f(n) + c \mbox{\ \ for all\ } n \in \nat.$$
\end{theorem}

\begin{proof} Note that, in order to make $A$ complex, we only need to complete
the construction from the previous theorem for infinitely many intervals.
It follows that every r.e.\ set $A$, that is array nonrecursive with respect
to $\{I_k\}_{k \in \nat}$, is also complex. In \cite[Theorem 2.5]{DJS90} it is
shown that every array nonrecursive degree contains such a set, i.e., it 
contains an r.e.\ complex set.

For the converse we use \cite[Theorem 4.1]{DJS90}. It states that if $A$ is r.e.\ and does not
have array nonrecursive degree, then for every total function $g \leq_T A$ there is a total 
recursive approximation $\overline{g}(x,s)$ such that $\lim_s \overline{g}(x,s) = g(x)$ and
$|\{ s : \overline{g}(x,s) \not= \overline{g}(x,s+1) \}| \leq x$, for all $x \in \nat$.
Actually, in \cite{DJS90} this is only stated for 0/1-valued $g$, but the proof provides
the more general version. 

Let $A$ be r.e.\ and not of array nonrecursive degree.
Assume we are given any total recursive, nondecreasing, unbounded function $f$.
Let $m(x) = 1 + \max\{ n : f(n) \leq x\}$; $m$ is total recursive.
Let $g(x) = \chi_A \restriction m(x)$. Since $g$ is recursive in $A$, there is
a total recursive approximation $\overline{g}(x,s)$ as above. 

How can we describe $\chi_A \restriction n$\ ? Given $n$ we compute $n' = \min\{ x : m(x) > n\}$.
Then we simulate $\overline{g}(n',s)$ until it outputs $g(n')$, which gives us $\chi_A \restriction n$.
In order to perform the simulation we only need to know the exact number $x \leq n'$
of mindchanges of $\overline{g}(n',s)$. Thus, $\chi_A \restriction n$ is specified by the pair
$\langle x, n' \rangle$ which can be encoded by a string of length $\log n + 2 \log(x+1) + O(1)$.
Since $m(x-1) \leq n$ we have $f(n) \geq x$, by the definition of $m$. Thus, we get
$$C(\chi_A \restriction n) \leq \log n + 2\log(x+1) + O(1) \leq \log n + f(n) + O(1).$$
This completes the proof of the theorem.
\end{proof}

Note that Theorem \ref{th:compdeg} entails the following curious gap phenomenon.
For every r.e.\ degree ${\bf a}$ there are only two cases: 
\begin{itemize}
\item[(1)] There is an r.e.\ set $A \in {\bf a}$ such that\\
$(\exists^{\infty} n)[C(\chi_A \restriction n) \geq 2 \log n - O(1)]$.
\item[(2)] There is no r.e.\ set $A \in {\bf a}$ and $\epsilon >0$ such that\\
$(\exists^{\infty} n)[C(\chi_A \restriction n) \geq (1 + \epsilon) \log n - O(1)]$.
\end{itemize}

\section{The Instance Complexity Conjecture fails}

In this section we determine the least possible 
instance complexity of nonrecursive r.e.\ sets.
Here it is convenient to take $A$ as a subset of $\{0,1\}^*$. 
Clearly, if $A$ is recursive then $ic(x:A)$ is bounded
by a constant for all $x$. The next result (another gap theorem) shows that, 
for infinitely many $x$, $ic(x:A)$ 
must be at least logarithmic in $C(x)$ if $A$ is nonrecursive\footnote{
This result was previously announced by Tromp \cite{Tr93}.}.

\begin{theorem} If $ic(x:A) \leq \log C(x) - 1$ for almost
all $x$, then $A$ is recursive. \label{th:gap}
\end{theorem}

\begin{proof} Let $P_k = \{0,1\}^{\leq k}$ and let ${\cal P}(P_k)$ denote the set
of all subsets of $P_k$. Uniformly in $k$ we enumerate a finite set $B_k \subseteq \{0,1\}^*$.

\smallskip
\noindent
{\it Step $0$}: Let $S_k = {\cal P}(P_k)$ and $B_k = \emptyset$.\\
{\it Step $n+1$}: Search via dovetailing for $I \in S_k$, $x \in \{0,1\}^*$ and $s \in \nat$
such that $U_s(p,x) = \perp$
for all $p \in I$. If such an $I$ is found, then enumerate $x$ into $B_k$, 
remove $I$ from $S_k$, and go to step $n+2$. $\Box$

\smallskip
\noindent
Note that $B_k$ is nonempty, since $I = \emptyset$ trivially satisfies the condition for all $x,s$. Also,
at most $|{\cal P}(P_k)| = 2^{|P_k|}$ elements are enumerated into $B_k$ and $B_k$ is uniformly r.e.
Thus, there is a partial recursive function $\psi : \{0,1\}^* \times \{0,1\}^* \rightarrow
\{0,1\}^*$ such that $\psi(\{0,1\}^{|P_k|},\lambda) = B_k$ for all $k$.
In particular, $C_{\psi}(x) \leq |P_k| = 2^{k+1}-1$ for all $x \in B_k$.
Choose a constant $c$ such that $C(x) \leq C_{\psi}(x) + c$ for all $x$.

Let $A \subseteq \{0,1\}^*$ be given and suppose that $ic(x:A) \leq \log C(x) - 1$
for almost all $x$. Then $ic(x:A) \leq \log(C_{\psi}(x) + c) - 1$, so
$ic(x:A) \leq \log C_{\psi}(x)$ for almost all $x$.
Since  $\lfloor \log C_{\psi}(x) \rfloor \leq \lfloor \log(2^{k+1}-1) \rfloor = k$ for all $x \in B_k$,
we can choose $k$ large enough such that for all $x \in B_k$ we have $ic(x:A) \leq k$.

Thus, for each $x \in B_k$ there is $p \in P_k$ such that $\lambda z.\ U(p,z)$ is a total $A$-consistent
function with $U(p,x) = \chi_A(x)$. Let $I_0 = \{ p \in P_k : \lambda z.\ U(p,z)$ is a total 
$A$-consistent function$\}$. This set is nonempty since $B_k$ is nonempty. 

Now consider the construction of $B_k$: Note that $I_0$ cannot be removed from $S_k$.
Otherwise there exists $x \in B_k$ such that $C_{\psi}(x) \leq 2^{k+1}-1$ and 
$U(p,x) = \perp$ for all $p \in I_0$, i.e., $ic(x:A) > k$, contradicting the choice
of $k$. Since $I_0$ is never removed from $S_k$, it follows from the construction of $B_k$
that for every $x$ there is $p \in I_0$ with $U(p,x) = \chi_A(x)$. Thus, if we amalgamate
all of the functions $U(p,-)$ with $p \in I_0$, we get a recursive characteristic function of $A$, i.e.,
$A$ is recursive.
\end{proof}

We prove that the lower bound of Theorem \ref{th:gap} is tight even for nonrecursive {\it r.e.} sets.
This refutes the Instance Complexity Conjecture of 
Orponen et al.\ \cite{Or90, O+94}, \cite[Exercise 7.41]{LV93}, stating that every
nonrecursive r.e.\ set has hard instances. In contrast, our result together
with Theorem \ref{th:gap} shows that the true threshold between the instance complexity
of recursive and nonrecursive sets is $\log C(x)$ instead of $C(x)$.

\begin{theorem} \label{th:icc1} There is a nonrecursive r.e.\ set $A$ 
and a constant $c$ such that\\
\mbox{\ } \hspace{4cm}$ic(x:A) \leq \log C(x) + c$ for all $x$.
\end{theorem}

\begin{proof} It suffices to construct a nonrecursive  r.e.\ set $A$ and a partial recursive
function $\psi$ such that $ic_{\psi}(x:A) \leq \log C(x) + 2$ for almost all $x$. In the following we write
$\psi_p$ for $\lambda z. \psi(p,z)$.

Let $E_k = \{ x : C(x) < 2^k-2 \}$ for $k \geq 1$. We want to establish that $ic_{\psi}(x:A) \leq k$
for all $x \in E_k$. Let $M_k = \{ p_{k,1},\ldots,p_{k,2^k-2}\}$ denote the set of the first $2^k-2$
strings of length $k$. The idea is that every $\psi_{p_{k,i}}$ is $A$-consistent and for each 
$x \in E_k$ there is $p \in M_k$ such that $\psi_p$ witnesses that $ic_{\psi}(x:A) \leq k$. 
There is, however, some difficulty to combine this with the requirement to make $A$ nonrecursive. 

The basic idea to satisfy the latter requirement is as follows: For each $e \geq 1$ we establish 
a unique diagonalization value $d_e$, then we wait until $d_e$ is enumerated into $W_e$, if this ever 
happens we enumerate $d_e$ into $A$. Here $\{ W_e \}_{e \in \nat}$ is the standard r.e.\ listing of all
r.e.\ sets of strings. Hence, this strategy makes sure that $\overline{A}$ is not r.e., so $A$ is a 
nonrecursive r.e.\ set. 

Suppose that $d_e$ appears in $E_k$ before it appears in $W_e$. If we define $\psi_{p_{k,i}}(d_e) = 0$
for some $i$, then, since $\psi_{p_{k,i}}$ should be $A$-consistent, we can no longer enumerate 
$d_e$ into $A$. This threatens our diagonalization strategy. On the other hand, we certainly
should make sure that $ic_{\psi}(d_e:A) \leq k$.

\noindent
This conflict is solved by a finite-injury priority argument:

If $e \geq k$ and we are forced to define $\psi_{p_{k,i}}(d_e) = 0$, then we assign a new much larger
value to $d_e$ and try to diagonalize at this new value. Note that $d_e$ is changed only finitely
often, because there are only finitely many values which may appear in $E_k$ for some $k \leq e$.
Thus, the value of $d_e$ eventually stabilizes and the $e$-th diagonalization strategy goes through
with this final value.

If $e < k$ then we do not use $\psi_{p_{k,i}}$ to ensure that $ic_{\psi}(d_e:A) \leq k$.
Thus, we define $\psi_{p_k,i}(d_e) = \perp$, which certainly maintains the $A$-consistency. 
Instead we will have two special programs $\tau_{e,1}, \tau_{e,2}$ 
of length $e$ (which are not in $M_e$;
this is the reason why we have left out two strings) to witness that $ic_{\psi}(d_e:A) \leq e < k$.
More precisely, if the final $d_e$-value is not enumerated into $A$, then $\psi_{\tau_{e,1}}$
will be the correct function. If the final $d_e$-value is enumerated into $A$, then
$\psi_{\tau_{e,1}}$ will not be $A$-consistent but $\psi_{\tau_{e,2}}$ is used as a back-up
function.

\noindent
It remains to explain how only $|M_k|$ many programs can take care of all of the elements in
$E_k$, which may be up to $2^{|M_k|}-1$ many. We show in an example how two programs $p_1, p_2$
can take care of $3 = 2^2-1$ elements (for simplicity, we drop the distinction between numbers and strings): 
At the beginning $\psi_{p_1}, \psi_{p_2}$ are undefined. 
Now in step $s_1$ the first element $x_1 < s_1$ appears. We let $\psi_{p_1}(x) = \chi_A(x)$ 
for all $x \leq s_1$. In the following steps $s$ we define $\psi_{p_1}(s) = \perp$ until the 
second element $x_2$ appears, say at step $s_2 > x_2$. If $x_2 \leq s_1$ we do nothing. If 
$x_2 > s_1$ then we define $\psi_{p_2}(x) = \chi_A(x)$ for all $x \leq s_2$ and in the following
steps $t$ we define $\psi_{p_2}(t) = \perp$. The point is that $\psi_{p_2}$ also takes care
of $x_1$, thus we suspend the definition of $\psi_{p_1}$ until a third element $x_3$ appears
at step $s_3 > x_3$. If $x_3 > s_2$ then we resume the definition of $\psi_{p_1}$ and let
$\psi_{p_1}(x) = \chi_A(x)$ for all $s_2 < x \leq s_3$. For arguments $t > s_3$ we define
both function equal to $\perp$. Note that now $p_1$ and $p_2$ together take care of 
$x_1,x_2,x_3$.

This idea is easily generalized: Let $succ(\sigma)$ denote the lexicographical successor
of $\sigma$, i.e., if $\sigma = b_1\ldots b_n \not= 1^n$ then $succ(\sigma) = 
0^{i-1}1b_{i+1}\ldots b_n$ where $i = \min\{ j : b_j = 0\}$. Then the programs
$p_{k,i} \in M_k$ with $succ^{(m)}(0^{|M_k|})(i) = 1$ take care, if exactly $m$ elements
are enumerated into $E_k$. (In the implementation below we count only  those elements which are
not $d_e$-values for some $e < k$.) Note that, since $m \leq 2^{|M_k|}-1$, $succ$ is never
applied to $1^{|M_k|}$.

We now turn to the detailed implementation. First we fix some additional notation 
and conventions. Let $\langle -,- \rangle$ denote a recursive pairing function 
which is increasing in its second argument. We assume that elements of $E_k$ are 
enumerated in steps such that in each step at most one new element is enumerated; also
if $x$ is enumerated in step $s$ then $l(x)<s$.
$W_{e,s}$ is the finite set of strings which are enumerated into $W_e$ in at most
$s$ steps of computation. 

In the construction the variables $e,i,j,k,n,s,t$ denote numbers, and $p,x,z$
denote strings. In addition, the following variables are used: $\psi_{p,s}$
the finite portion of $\psi_p$ constructed up to stage $s$; the $i$-th bit of
$\sigma_{k,s} \in \{0,1\}^{|M_k|}$ tells us whether $\psi_{p_{k,i}}$ is currently 
assigned to take care of the elements in $E_k$; $len(k,s)$ is the greatest length $n$
such that our set-up at stage $s$ guarantees that $ic_{\psi}(x:A) \leq \log C^s(x)+2$
for all $x \in E_{k,s}$ with $l(x)< n$; $d_e(s)$ is the current value of the $e$-th diagonalization
point. We call $e$ ``active'' as long as no $d_e$-value has been enumerated into $A$,
otherwise we call $e$ ``passive''. So, if $e$ is ``passive'', then we know that we have
explicitly satisfied the $e$-th diagonalization requirement. $A_s$ denotes the finite 
set of elements which have been enumerated into $A$ up to stage $s$.

Let $R(k,s) = \{ d_e(s') : e < k\ \wedge\ s' \leq s \}$. As explained above, the programs
in $M_k$ do not need to take care of the elements in $R(k,s)$.

If one of the variables $v(s)$ is not explicitly changed at stage $s+1$, then we
assume without further mentioning that $v(s+1) = v(s)$.

We first describe the construction of $\psi_{p}$ for $p \in M_k$, $k \geq 1$. Then
we define $\psi_p$ for the two special values $p = \tau_{e,1}, \tau_{e,2}$ of each length $e$.
\begin{tabbing}
{\it Construction:}\\
{\it Stage $0$}: \= Let $\psi_{p} = \lambda x. \uparrow$, for all $p \in \{0,1\}^*$.
For all $k \geq 1$: $\sigma_{k,0} = 0^{|M_k|}; len(k,0) = 0$;\\
\> $d_k(0) = 0^{\langle k,0 \rangle}$;
declare $k$ as ``active''. Let $A_0 = \emptyset$.\\
{\it Stage $s+1$}:\\
Case I: $s$ is even.\\
\quad \= For $e = 0,\ldots,s$: If $e$ is active and $d_e(s) \in W_{e,s} - A_s$, then enumerate
$d_e(s)$ into $A$\\
\> and declare $e$ ``passive''.\\
Case II: $s$ is odd, $s = 2 \langle k,t \rangle + 1$.\\
\> Let $\psi_{p_{k,i}}(x) = \perp$, for all $i$ with $\sigma_{k,s}(i) =1$ and all $x$ with $l(x) = t$.\\ 
\> If a new element $x$, $l(x) < t$, enters $E_k$ after exactly $t$ steps, then act\\
\>  according to the following cases:\\
\quad\quad\= a.) If $x \in R(k,s)$ or $l(x) < len(k,s)$ then go to stage $s+2$.\\
\> b.) \= Otherwise do the following:\\ 
\> \>Let $\sigma_{k,s+1} = succ(\sigma_{k,s})$ and  $i = \min\{j : \sigma_{k,s+1}(j) = 1\}$. 
   (At at most $2^{|M_k|}-1$\\
\> \>  elements are enumerated in $E_k$, so we get $\sigma_{k,s+1} \in \{0,1\}^{|M_k|}-\{0^{|M_k|}\}$.)\\ 
\> \> Let $n = \min\{ l(z) : z \not\in dom(\psi_{p_{k,i}})\}$.\\
\> \> Define $\psi_{p_{k,i}}(z) = \chi_{A_s}(z)$, for all $z \not\in R(k,s)$ such that $n \leq l(z) \leq t$.\\
\> \> Let $\psi_{p_{k,i}}(z) = \perp$, for all $z \in R(k,s)$ such that $n \leq l(z) \leq t$.\\
\> \> Let $len(k,s+1) = t+1$. For all active $e \geq k$, let $d_e(s+1) = 0^{\langle e,s+1 \rangle}$.\\
\> \> Go to stage $s+2$. $\Box$
\end{tabbing}

\noindent
For each $e \geq 1$ we define

$$\psi_{\tau_{e,1}}(x) = \cases{0         &if $x \in range(d_e)$;\cr
                               \perp   &otherwise.\cr}$$

If $e$ is active at all stages then let $\psi_{\tau_{e,2}} = \lambda x. \uparrow$.
Otherwise let $s_e$ be the (unique) stage where $e$ is declared ``passive'' and let

$$\psi_{\tau_{e,2}}(x) = \cases{0         &if $(\exists t < s_e)[x = d_e(t) \not= d_e(s_e)]$;\cr
                                1         &$x = d_e(s_e)$;\cr
                               \perp   & otherwise.}$$
{\it End of Construction.} 

\medskip
\noindent
{\it Verification:}\\
Most of the following claims are standard. The crucial one is Claim 3, b.), c.).

\begin{tabbing}
{\it Claim 1\/} For all $e \geq 1$:\\
a.) \=  $l(d_e)$ is nondecreasing, and for all $s$: If $d_e(s) \not= d_e(s+1)$ then $l(d_e(s+1)) > s$.\\
b.)\> $range(d_e)$ is a uniformly recursive finite set; $range(d_e) \cap range(d_{e'}) = \emptyset$\\ 
\> for all $e' \not= e$.\\
c.)\> If $A \cap range(d_e)$ contains an element $x$ then $x = \lim_{s \rightarrow \infty} d_e(s)$.\\
d.)\> For all $x,s$: If $l(x) \leq s$ and, for all $e$,  $x \not= d_e(s)$, then $x \not= d_e(s')$
for all $e$\\
\> and all $s' \geq s$.\\
e.)\> $A$ is r.e.\ and nonrecursive.
\end{tabbing}

\noindent
{\it Proof:\/} a.) If $d_e(s) \not= d_e(s+1)$ then for some $s' \leq s$:
$l(d_e(s)) = \langle e,s' \rangle < \langle e,s+1 \rangle = l(d_e(s+1))$.
Note that $\langle e,s+1 \rangle > s$, since $\langle -,- \rangle$ is monotone in the
second argument.

b.) It follows from a.) that $range(d_e)$ is uniformly recursive. It is a finite set,
because $d_e(s)$ changes only if a new element is enumerated in some set $E_k$, $k<e$
which happens only finitely often. So $\lim_{s \rightarrow \infty} d_e(s)$ exists and is finite.
$range(d_e)$ and $range(d_{e'})$ are disjoint for $e \not= e'$, since $\langle -,- \rangle$
is injective.

c.) If $d_e(s)$ is enumerated into $A$ at stage $s+1$, then $e$ is declared ``passive'', so $d_e(s)$ is
fixed at all later stages.

d.) This follows from a.).

e.) Clearly $A$ is r.e. Suppose for a contradiction that $A$ is recursive.
Then there exists $e$ with $\overline{A} = W_e$. By a.), b.), there is a stage $s$ such that
$d_e(s') = d_e(s)$ for all $s' \geq s$. By construction, $d_e(s)$ is enumerated into $A$
iff it is enumerated into $W_e$. This contradicts the hypothesis $\overline{A} = W_e$. $\Box$

\begin{tabbing}
{\it Claim 2\/} For all $e \geq 1$:\\
a.) \=  $\psi_{\tau_{e,1}}, \psi_{\tau_{e,2}}$ are uniformly partial recursive.\\
b.)\> If $e$ is always ``active'' then $\psi_{\tau_{e,1}}$ witnesses that $ic_{\psi}(x:A) \leq e$ for all\\
\> $x \in range(d_e)$.\\
c.)\> If $e$ is eventually ``passive'' then $\psi_{\tau_{e,2}}$ witnesses that $ic_{\psi}(x:A) \leq e$ for all\\
\> $x \in range(d_e)$.
\end{tabbing}

\noindent
{\it Proof:\/} a.) This follows from Claim 1, b.)

b.) If $e$ is always ``active'' then $range(d_e) \cap A = \emptyset$, thus
$\psi_{\tau_{e,1}}$ is $A$-consistent and $\psi_{\tau_{e,1}} = 0 = \chi_A(x)$ 
for all $x \in range(d_e)$.

c.) If $e$ is declared ``passive'' at stage $s+1$ then $A \cap range(d_e) = \{ d_e(s) \}$.
Thus $\psi_{\tau_{e,2}}$ is $A$-consistent and $\psi_{\tau_{e,2}}(x) = \chi_A(x)$ for all
$x \in range(d_e)$. $\Box$

\medskip
\noindent
Let $\psi^s_e$ denote the finite portion of $\psi_e$ defined at the end of stage $s$.

\begin{tabbing}
{\it Claim 3\/} For all $s = 2\langle k,t \rangle +1$:\\
a.) \= For all $i$, $1 \leq i \leq |M_k|$: $\psi^s_{p_{k,i}}$ is an
$A$-consistent function.\\
b.) \> For all $i$, $1 \leq i \leq |M_k|$: If $\sigma_{k,s+1}(i) = 1$ then
$dom(\psi^{s+1}_{p_{k,i}}) = \{ x : l(x) \leq t \}$.\\
c.) \> For all  $x$, $l(x) < len(k,s+1)$: If $x \not \in R(k,s)$ then there exists
$i$, $1 \leq i \leq |M_k|$,\\
\> with $\sigma_{k,s+1}(i) = 1$ and $\psi^{s+1}_{p_{k,i}}(x) = \chi_{A_{s+1}}(x)$.
\end{tabbing}

\noindent
{\it Proof:\/} a.) We use  Claim 1, d.) and the fact that $\psi_{p_{k,i}}$ is defined at stage
$s+1$ only for arguments less than $s$. If $e<k$ and $\psi_{p_{k,i}}(d_e(s'))$ is defined
then $\psi_{p_{k,i}}(d_e(s')) = \perp$, so there is no problem with consistency. If $e \geq k$
and $\psi_{p_{k,i}}(d_e(s')) = \chi_{A_{s+1}}(d_e(s'))$ is defined at stage $s+1 > s'$, then
either $e$ is already ``passive'', so $\chi_{A_{s+1}}(d_e(s')) = \chi_A(d_e(s'))$, or $e$ is ``active''
and we define $d_e(s+1)$ at stage $s+1$ such that $l(d_e(s+1)) > l(d_e(s'))$. In the latter
case we get $\psi_{p_{k,i}}(d_e(s')) = \chi_{A_{s+1}}(d_e(s')) = \chi_A(d_e(s')) = 0$.

\smallskip
\noindent
b.) and c.) are shown by induction on $s$. Consider stage $s+1 = 2 \langle k,t \rangle + 2$.
If no new element is enumerated in $E_k$ after exactly $t$ steps then $\sigma_{k,s+1} = \sigma_{k,s}$
and b.), c.) follow from the induction hypothesis and the definition of $\psi_{p_{k,i}}$ at stage $s+1$.

Now assume that $x$ enters $E_k$ after exactly $t$ steps.
If case a.) occurs, the claim follows from the induction hypothesis. 
If case b.) occurs, we have $x \not\in R(k,s)$ and $len(k,s) \leq l(x) < t$.
We have $\sigma_{k,s+1} = succ(\sigma_{k,s})$, so $\sigma_{k,s+1}(i') = \sigma_{k,s}(i')$ for all
$i' > i = \min\{ j : \sigma_{k,s+1}(j) = 1\}$. 

If $\sigma_{k,s+1}(i') = 0$ for all $i' > i$ then
$s+1$ is the first stage $s'$ where $\sigma_{k,s'}(i) = 1$.
This means that $\psi^s_{p_{k,i}} = \lambda x. \uparrow$ and 
$\psi^{s+1}_{p_{k,i}}(z) = \chi_{A_{s}}(z) = \chi_{A_{s+1}}(z)$,
for all $z$ such that $l(z) \leq t$ and $z \not\in R(k,s)$.

If there is $i' > i$ with $\sigma_{k,s+1}(i') = 1$, then there exists a greatest stage $s' < s$ with
$\sigma_{k,s'}(i) = 1\ \wedge \sigma_{k,s'+1}(i) = 0$. Then we have $\sigma_{k,s'+1}(i') =
\sigma_{k,s+1}(i')$ for all $i'>i$ and we have $\sigma_{k,s'+1}(i') = \sigma_{k,s+1}(i') = 0$ for all $i' < i$.
By induction hypothesis, we get for all $x$ with $l(x) < len(k,s'+1)$:
If $x \not\in R(k,s')$ then there exists $j$ with $\sigma_{k,s'+1}(j) = 1$ and $\psi^{s'+1}_{p_{k,j}}(x)
= \chi_{A_{s'+1}}(x) = \chi_A(x)$ (the second equality holds by part a.)). 

Since $R(k,s') \subseteq R(k,s)$,  it only remains to consider $x$
with $len(k,s'+1) \leq l(x) < len(k,s+1) = t+1$. As $\sigma_{k,s'}(i) = 1$ it follows, by induction
hypothesis, that $dom(\psi^s_{p_{k,i}}) = dom (\psi^{s'}_{p_{k,i}}) = \{ x : l(x) < len(k,s'+1)\}$.
Thus, $n = \min\{ l(z) : z \not\in dom(\psi^s_{p_{k,i}})\} = len(k,s'+1)$ and at stage $s+1$ we define
$\psi^{s+1}_{p_{k,i}}(z) = \chi_A(z)$, for all $z \not\in R(k,s)$ such that 
$len(k,s'+1) = n \leq l(z) < t+1 = len(k,s+1)$. 
For $z \in R(k,s)$ and $l(z) \leq t$ we have $\psi^{s+1}_{p_{k,i}}(z) = \perp$.
This completes the proof of b.), c.). $\Box$

\medskip
\noindent
{\it Claim 4\/} For almost all $x$: $ic_{\psi}(x:A) \leq \log C(x) + 2$. 

\smallskip
\noindent
{\it Proof:\/} Let $k \geq 1$ be minimal such that $x \in E_k$. If $x \in R(k,s)$ for some $s$ then,
by Claim~2, we get $ic_{\psi}(x:A) < k$. If $x \not\in R(k,s)$ for all $s$ then let 
$\sigma_k = \lim_{s \rightarrow \infty} \sigma_{k,s}$. By Claim 3, there exists $i, 1 \leq i \leq |M_k|$ 
such that $\sigma_k(i) = 1\ \wedge\ \psi_{p_{k,i}}(x) = \chi_A(x)$. Furthermore, $\psi_{p_{k,i}}$ is
total recursive and $A$-consistent, so $ic_{\psi}(x:A) \leq k$. Since $E_1 = \emptyset$, we have $k>1$
and $x \not\in E_{k-1}$, so $2^{k-1}-2 \leq C(x)$, i.e., $k \leq \log(C(x)+2)+1 \leq \log C(x) +2$
for all $x$ with $C(x) \geq 2$. $\Box$
\end{proof}

\noindent
What happens for $\overline{ic}$ ?
Of course, the instance complexity conjecture also fails for $\overline{ic}$.
It even fails in a much stronger way, because, in contrast to  Theorem \ref{th:gap}, 
$\overline{ic}$ can be arbitrary small, as we now show.

\begin{theorem} \label{th:icc2} For every recursive function $f$ there is an r.e.\ nonrecursive set $A$
such that $$f(\overline{ic}(x:A)) \leq C(x) \mbox{\ for almost all\ } x.$$
\end{theorem}

\begin{sketch}
We may assume that $f$ is strictly increasing. As above it suffices to 
define a partial recursive function $\psi(p,x)$ such that 
$f(\overline{ic}_{\psi}(x:A)) \leq C(x)$ for almost all $x$ and $A$ is nonrecursive. 
This leads to the following requirements for all $i \geq 1$:
\begin{tabbing}
$(N_i)$\ \= $(\forall x) [ C(x) < f(i+1)\ \Rightarrow\ (\exists p \in \{0,1\}^i)
[ \chi_A \mbox{\ extends\ } \psi_p \mbox{\ and\ } \psi_p(x) = \chi_A(x) ]].$\\
$(P_i)$ \> $W_i \not= \overline{A}.$
\end{tabbing}
These can be satisfied by an easy finite-injury construction.
Fix an enumeration of $E_i = \{ x : C(x) < f(i+1) \}$ for all $i$.

During the construction we have for each $i$ a current $p_i \in \{0,1\}^i$ 
which satisfies $(N_i)$ for all $x$ that have been currently enumerated into $A$.
If some $x$ with $\psi_{p_i}(x) = 0$ is later enumerated into $A$, then $\psi_{p_i}$
is no longer $A$-consistent and we have to choose a new $p_i$. Since we have $2^i$
candidates for $p_i$, we can afford $2^i-1$ injuries.

Therefore, we allow to enumerate a diagonalization witness $x$ into $A$ at stage $s$
for the sake of $(P_i)$, only if $x$ has not yet appeared in any $E_j$ with $j \leq i$.
Clearly, $(P_i)$ can still be satisfied. Furthermore, $(N_i)$ is injured at most
$i$ times. Since $i \leq 2^i-1$ for all $i\geq 1$, every $(N_i)$ will be eventually satisfied.
\end{sketch}

\noindent
{\it Remark:\/} In the course of the construction at most $2^{f(i+1)}-1$ elements 
are not allowed to be enumerated into $A$ by $(P_i)$. Hence, we can fix in advance
a set $J_i$ of $2^{f(i+1)}$ witnesses for $(P_i)$ and guarantee that one of them will
be successful. Therefore, we can also modify the construction and satisfy the
following requirements $(P'_i)$ instead of $(P_i)$ for any fixed r.e.\ set $B$
$$(P'_i)\ \ i \in B \ \Leftrightarrow\  J_i \cap A \not= \emptyset.$$
Then we get $B \leq_d A$. If we choose $B = K$, this shows that there
is a d-complete set which satisfies the condition of the theorem.
Since we need to enumerate at most one element of $J_i$ into $A$,
we get that $A \leq_{wtt(1)} B$. Thus, every r.e.\ wtt-degree contains
a set $A$ as in the theorem. It can be shown 
that this does not hold for r.e.\ tt-degrees.

\section{R.e.\ sets having hard instances}

While we have shown in the last section that ICC fails for some nonrecursive r.e.\ sets,
it is interesting to find out whether there are properties of r.e.\ sets which imply the 
existence of hard instances. We consider this question for classes of complete sets and of
simple sets. Indeed, in most cases it turns out that such sets must have hard instances,
which is a partial resurrection of ICC.

Buhrman and Orponen \cite{BO93}, \cite[Exercise 7.40]{LV93} proved that
the set of all random strings $R = \{ x : C(x) \geq l(x) \}$ satisfies 
$ic(x:R) \geq l(x) - O(1)$ for all  $x \in R$. (Actually, their result also 
holds for $\overline{ic}$ instead of $ic$.) 
Using the observation  
$$\mbox{(*)}\quad \mbox{If\ }A \leq_m B \mbox{\ via\ } f \mbox{, then\ }
ic(x:A) \leq ic(f(x):B) + O(1) \mbox{ \ for all\ } x.$$
and the fact that $R$ is co-r.e., they conclude that every 
m-complete set $A$ has hard instances in its complement. 
They asked whether the hard instances can be chosen 
from $A$ instead of $\overline{A}$. 
(This is of course impossible in the $\overline{ic}$-version.)
The next result gives a positive answer.

\begin{theorem} There is an r.e.\ set $A$ with $ic(x:A) \geq l(x)$ 
for infinitely many $x \in A$. 
\end{theorem}

\begin{proof} Uniformly in $n$ we enumerate $A \cap \{0,1\}^n$ as follows:
Let $x_1,\ldots,x_{2^n}$ be a listing of all strings of
length $n$ in lexicographical order.

\begin{tabbing}
{\it Step $0$}: Enumerate $x_1$ into $A$, let $i = 1$, $I = \{0,1\}^{\leq n-1}$,
$J = \{1,\ldots,2^n\} - \{1\}$.\\
{\it Step $s+1$}: If there is a string $p \in I$ such that\\
\quad\quad \= (a) $U_s(p,x_j) \in \{0,1\}$ for some $j \in J$, or\\
\>       (b) $U_s(p,x_j) = \perp$ for all $j \in J$,\\
then choose the least such $p$, let $I = I -\{p\}$, and do the following:\\
\> In case (a): Enumerate $x_j$ into $A$ iff $U_s(p,x_j) = 0$. Let $J = J - \{j\}$.\\
\> In case (b): Let $i = \min(J)$. Enumerate $x_i$ into $A$ and let $J = J - \{i\}$. $\Box$
\end{tabbing}

\noindent
At the end of Step 0 we have $|I| = |J| = 2^n-1$. In all later steps an element of
$I$ is removed iff an element of $J$ is removed. Thus, at the end of each step
we have $|I| = |J|$. Also, if case (b) occurs then $\min(J)$ exists (since at that point
$|J| > 0$). Note that the value of $\chi_A(x_j)$ is fixed when $j$ is removed from $J$.

Let $i_0, I_0, J_0$ be the final values of $i, I, J$ in the above construction
and choose $s_0$ such that $i = i_0, I = I_0, J = J_0$ in all steps $t \geq s_0$.
Suppose for a contradiction that $ic(x_{i_0}:A) < n$ via $p \in \{0,1\}^{\leq n-1}$.

If $p \not\in I_0$ then there is a stage $s \leq s_0$ when $p$ was removed from $I$.
If $p$ was removed in case (a) via $j$, then $U(p,x_j) \not= \chi_A(x_j)$.
If $p$ was removed in case (b) then $U(p,x_{i_0}) = \perp$. Hence, $p$ 
does not witness that $ic(x_{i_0}:A) < n$, a contradiction.

If $p \in I_0$ then  $|J_0| = |I_0| \geq 1$ and there is $t>s_0$ such that
$U_t(p,x) \in \{0,1,\perp\}$ for all $x \in J_0$. Hence, at stage $t+1$ either case (a)
or case (b) occurs and $|I_0|$ decreases, contradicting the choice of $s_0$.

Thus, we have $ic(x_{i_0}:A) \geq n = l(x_{i_0})$ and clearly $x_{i_0} \in A$. 
Since this holds for all $n$, the theorem is proved.
\end{proof}

Using (*) we get the following corollary.

\begin{corollary} For every m-complete set $A$ 
there is a constant $c$ such that\\
\mbox{\ } \hspace{3.5cm}$ic(x:A) \geq C(x) - c$ for infinitely many $x \in A$.
\end{corollary}

This result also holds for a much weaker reducibility, as we now show.

\begin{theorem} For every  wtt-complete set $A$ 
there is a constant $c$ such that\\
\mbox{\ } \hspace{3.5cm}$ic(x:A) \geq C(x) - c$ for infinitely many $x \in A$.
\end{theorem}

\begin{proof} Suppose that $A$ is a wtt-complete set. We enumerate an auxiliary
r.e.\ set $B$ and a uniformly r.e.\ sequence $\{E_n\}_{n \in \nat}$ with $|E_n| \leq 2^n$.
Then there is a partial recursive function 
$\psi: \{0,1\}^*\times \{0,1\}^* \rightarrow \nat$ such that
$\psi(\{0,1\}^n, \lambda) = E_n$. Hence, $C_{\psi}(x) \leq n$ for all $x \in E_n$ and there is a
constant $c$, independent of $n$, such that $C(x) \leq n+c$ for all $x \in E_n$.
Thus, it suffices to satisfy the following requirement for all $n$
$$(R_n)\ \ (\exists x \in E_n \cap A)[ ic(x:A) \geq n-1].$$
By the recursion theorem and the fact that $A$ is wtt-complete, we can assume that we
are given in advance the index of a wtt-reduction from $B$ to $A$, i.e., 
a Turing reduction $\Phi$ and a total recursive use-bound $g$ such that, for all $x$,
$\chi_B(x) = \Phi^A(x)$ and in the computation of $\Phi^A(x)$ every query is less
than $g(x)$. 

Each $(R_n)$ is satisfied independently from the other requirements;
so for the following fix $n$ and let 
$x_1 = \langle n, 1 \rangle, \ldots, x_{2^n} = \langle n, 2^n\rangle$,
$m = \max\{ g(x_i) : 1 \leq i \leq 2^n\}$, and  $I = \{ p : l(p) < n-1\}$.
We enumerate $E_n$ and $B \cap \{x_1,\ldots,x_{2^n}\}$ in steps $i = 0,\ldots, 2^n$
as follows:

\medskip
\noindent
{\it Step $0$}: Let $s_0 = 0, E_n = \emptyset$.\\
{\it Step $i+1$}: Search for the least $s \geq s_i$ such that
\begin{enumerate}
\item[(1)] $\Phi^{A_s}_s(x_j) = 1$ with use less than $g(x_j)$ for $j = 1,\ldots,i$ and
$\Phi^{A_s}_s(x_j) = 0$ with use less than $g(x_j)$ for $j = i+1,\ldots,2^n$.
\item[(2)] For each $x \in E_n$ there is $p \in I$ such that
\begin{enumerate}
\item[(2.1)] $U_s(p,z)$ is defined for all $z \leq m$.
\item[(2.2)] $(\forall z \leq m)[ U_s(p,z) \not= \perp\ \Rightarrow\ U_s(p,z) = \chi_{A_s}(z)].$
\item[(2.3)] $U_s(p,x) = 1$.
\end{enumerate}
\end{enumerate}
Let $s_{i+1} = s$. 
Enumerate $x_{i+1}$ into $B$ and compute some $x \leq m$ with $x \in A-A_{s_{i+1}}$.
(Note that $x$ exists because otherwise $\Phi^{A}(x_{i+1}) = 0 \not= 1 = \chi_B(x_{i+1})$.
We can find $x$ by enumerating $A$.) Let CONS be the set of all $p \in I$ 
which satisfy conditions (2.1) and (2.2) for $s = s_{i+1}$.
If  $U_s(p,x) = \perp$ for all $p \in$ CONS, then enumerate $x$ into $E_n$.
Goto step $i+2$.
$\Box$
 
\medskip
\noindent
By construction, we have $E_n \subseteq A$.
We want to argue that in some step of the construction the search does not terminate. 
Since $\chi_B(x) = \Phi^A(x)$, this can only happen if condition (2) is
not satisfied for any sufficiently large $s$. But this means that $ic(x:A) \geq n-1$
for some $x \in E_n$.

Consider the value of CONS $\subseteq I$ after each terminating step: 
We show that a new element enters CONS or an element is removed forever from CONS.
Since there are at most $|I| < 2^{n-1}$ strings which may at some point
become a member of CONS, it follows that there are less than $2 \cdot 2^{n-1} = 2^n$
terminating steps, which completes the proof.

Note that if a string $p$ is removed from CONS at some stage $s$,
then there is $x$ such that $U_s(p,x) = 0$ and $\chi_{A_s}(x) = 1$.
Thus, $x$ cannot enter CONS again at any later stage.

Suppose that step $i+1$ terminates and consider the current value of CONS and of $x$
at the end of this step. There are two cases:

(a) $U_s(p,x) = \perp$ for all $p \in$ CONS. Then $x$ is enumerated into $E_n$,
so in the next step a new string must enter CONS such that condition (2.3) is satisfied
for $x$.

(b) $U_s(p,x) \not= \perp$ for all $p \in$ CONS. Hence, $U_s(p,x) = 0$ and, since
$\chi_{A_{s_{i+2}}}(x) = 1$, $p$ is removed from CONS if the next step terminates.
\end{proof}

By a similar proof, one can show that every btt-complete set has hard instances
w.r.t.\ $\overline{ic}$. We have noticed
in the remark following Theorem \ref{th:icc2} that this 
is no longer true for d-complete sets. 
But we can show that it still holds for Q-complete sets.

Recall that $A$ is Q-complete if it is r.e.\ and there is a recursive function $g$
such that for all $x$:
$$x \in K\ \Leftrightarrow\ W_{g(x)} \subseteq A.$$
See \cite[p.\ 281 f.]{Od89} for more information on Q-reducibility.

\begin{theorem} Every Q-complete set $A$ has hard instances, even w.r.t.\ $\overline{ic}$.
\end{theorem}

\begin{proof}  Suppose that $A$ is Q-complete. As in the previous proof we enumerate an auxiliary
r.e.\ set $B$ and an r.e.\ sequence of finite sets $\{E_n\}_{n \in \nat}$ such that $|E_n| \leq 2^n$.
It suffices to get infinitely many $n$ such that there is $y \in E_n$ with
$\overline{ic}(y:A) \geq n-2$.

By the recursion theorem and the Q-completeness of $A$, we may assume that we are given in advance
a recursive function $g$ such that $B \leq_Q A$ via $g$, i.e., for all $x$, $x \in B \Leftrightarrow 
W_{g(x)} \subseteq A$. 

The first idea is to run a version of the previous construction: 
We keep a number $x$ out of $B$ and find $y \in W_{g(x)}$ which has not yet been enumerated 
into $A$. Then we enumerate $y$ into $E_n$ and wait until some $A$-consistent program $p$
with $l(p) < n-2$ shows up and $U(p,y) = 0$. Then we enumerate $x$ into $B$,  which forces
$y$ into $A$ and diagonalizes $p$.

However, this approach does not work, because it might happen that after we enumerate
$y$ into $E_n$, $y$ is also enumerated into $A$, and {\it after that\/} $U(p,y) = 1$ is defined. 
Then we cannot diagonalize $p$ by enumerating $x$ into $B$, but we have incremented $|E_n|$. 
Since this can happen an arbitrary finite number of times, we run into conflict with 
the requirement $|E_n| \leq 2^n$.

Therefore, we use the following modification: For each $n$, if $E_n \not=\emptyset$ then
we enumerate $y$ into $E_n$ only if $y$ has been previously enumerated into $E_{n+1}$,
and then we proceed according to the first idea. If later $y$ is enumerated into $A$ we get
a diagonalization for $n+1$ instead of $n$, which is also fine.

Now we turn to the formal details: Let $I_n = \{ p : l(p) < n-2 \}$.
$p \in \{0,1\}^*$ is called {\it $A$-consistent\/} at stage $s+1$ if, for all $z \leq s$,
either $U_s(p,z)$ is undefined or $U_s(p,z) = \chi_{A_s}(z)$.

\noindent
We maintain the following invariant for all $n,s,y$:\\
If $E_n \not=\emptyset$ at stage $s+1$ then enumerate $y$ into $E_n$ only if $P(n,s,y)$ holds,
where:

\begin{tabbing}
\quad\quad $P(n,s,y)\ \Leftrightarrow$\ \= $y \in E_{n+1}-A_s$, $E_{n} \subseteq A_s$, and there 
is $p \in I_{n+1}$ which\\
\> is $A$-consistent at stage $s+1$ and $U_s(p,y) = 0$.
\end{tabbing}

As a consequence of this invariant it already follows that $|E_n| \leq 2^n$:
Suppose that $E_n \not= \emptyset$ and we enumerate $y$ into $E_n$ at stage $s+1$.
Then we enumerate the next element into $E_n$ only after $y$ has been enumerated 
into $A$, and hence the program $p \in I_{n+1}$ which had witnessed the condition 
$P(n,s,y)$ is diagonalized and can never be $A$-consistent again. 
Since $|I_{n+1}|< 2^{n-1}$, it follows that we will enumerate at 
most $1+2^{n-1}$ programs into $E_n$. In particular, $|E_n| \leq 2^n$ for all $n$.

We say that $n$ is {\it saturated\/} at stage $s+1$ if, for every $y \in E_n$, there is $p \in I_n$
such that $p$ is $A$-consistent at stage $s+1$ and $U_s(p,y) = \chi_{A_s}(y)$.
The goal of the construction is to produce infinitely many $n$ which are almost
always not saturated. This implies at once that there are infinitely many $y \in E_n$
with $\overline{ic}(y:A) \geq n-2$, and we are done.

To achieve this goal we construct a sequence $d_0 < d_1 < d_2 < \cdots$ and satisfy the
following requirements
$$(R_i)\ \ \mbox{The interval\ } [d_i,d_{i+1}) \mbox{\ contains an\ } n 
\mbox{\ which is almost always not saturated}.$$
The $d_i$'s are constructed by recursive approximation: The value of $d_i$ may change
finitely often and eventually stabilizes. Some additional variables are needed for
book-keeping: For each $i$ there is a finite set $T_i$ containing the set of all $x$
which may be enumerated into $B$ for the sake of $(R_i)$. For each $n$ we have
three variables $active(n), cand(n), source(n)$. $active(n)$ is a Boolean flag
which indicates if there is some $y \in E_{n}-A_s$ to be enumerated into $E_{n-1}$;
in this case $cand(n) = y$ and $source(n) = x$ such that $x \not\in B_s$ 
and $y \in W_{g(x),s}$.

We say that $i$ {\it requires attention\/} at stage $s+1$ if one of the following conditions holds
at the beginning of of stage $s+1$.
\begin{enumerate}
\item[(1)] $d_{i+1}$ is undefined.
\item[(2)] $d_{i+1}$ is defined and every $n \in [d_i,d_{i+1})$ is saturated at stage $s+1$.
\end{enumerate}

\noindent
{\it Construction:}\\
{\it Stage $0$}: Let $d_0 = 0, d_{i+1} = \uparrow, T_i = \emptyset$ for all $i$. Let $active(n) = 0, 
E_n = \emptyset$ for all $n$.\\
{\it Stage $s+1$}: For every $n$ such that $active(n) = 1$ and $cand(n) \in A_s$ let $active(n) = 0$.\\
Let $i$ be the least number which requires attention at stage $s+1$. 
If it requires attention through (1) then let $d_{i+1} = s+1$.\\ 
If it requires attention through (2) then we distinguish two cases:

(a) If there is a least $n \in (d_i, d_{i+1})$ such that $active(n) = 1$ and $E_{n-1} \subseteq A_s$,
then enumerate $cand(n)$ into $E_{n-1}$ and let $active(n) = 0$. 
If $n-1 = d_i$ then enumerate $source(n)$ into $B$, else let $active(n-1) = 1$,
$cand(n-1) = cand(n)$, and $source(n-1) = source(n)$.

(b) Otherwise put $s+1$ into $T_i$ and let $x = \min(T_i-B_s)$. Find the least $s'$ such that
$W_{g(x),s'}-A_s \not=\emptyset$ and let $y = \min(W_{g(x),s'}-A_s)$. 
Let $active(s+1) = 1, cand(s+1) = y, source(s+1) = x$, and enumerate $y$ into $E_{s+1}$.

In both cases let $T_i = T_i \cup \bigcup_{j>i} T_j$ 
and let $T_j = \emptyset, d_j = \uparrow$, for all $j>i$.\\
{\it End of Construction.}

\medskip
\noindent
It easily follows by induction on $s$ that our invariant is satisfied: Note that before we enumerate
a new number into $E_{n-1}$ via step (a), we require that $E_{n-1} \subseteq A_s$. If we enumerate
a number via step (b) then the corresponding set was previously empty. Therefore, at each stage $s+1$
every $E_n$ contains at most one number which is not in $A_s$. Now suppose that $E_{n-1} \not=\emptyset$
at the end of stage $s$ and we enumerate a number $y$ into $E_{n-1}$ at stage $s+1$. Then case (a)
occurred and $y = cand(n) \not\in A_s$ (since $active(n) = 1$). By the previous remarks, we have
$E_{n-1}\subseteq A_s$. Since $n$ is saturated at stage $s+1$, there is an $A$-consistent $p \in I_n$
such that $U_s(p,y) = \chi_{A_s}(y) = 0$. Thus, $P(n-1,s,y)$ holds.

Hence, it only remains to verify that requirement $(R_i)$ is satisfied for all $i$. This is done by induction
on $i$. By induction hypothesis, there is a least stage $s_0$ such that $d_i = s_0$ 
is defined at stage $s_0$ and no $i'<i$ requires attention at any stage $s > s_0$. 
At the end of stage $s_0$ we have  $E_{d_i} = \emptyset$ and $T_i = \emptyset$.
We have shown above that the cardinality of $E_{d_i}$ is always bounded by $2^{d_i}$. Hence, there
exists $s_1 \geq s_0$  such that $E_{d_i}$ does not change after stage $s_1$. 
Note that $E_{d_i} \subseteq A$, because each time when we enumerate $y$ into $E_{d_i}$,
we enumerate some $x$ into $B$ such that $x\in B \Leftrightarrow W_{g(x)}\subseteq A$
and $y \in W_{g(x)}$; thus we force $y$ into $A$. So we can choose $s_1$ large enough
such that $E_{d_i} \subseteq A_{s}$ for all $s \geq s_1$.

Suppose for a contradiction that $i$ requires attention infinitely often. We will argue that at some
stage $s_2 > s_1$ a new element is enumerated into $E_{d_i}$, which contradicts the choice of $s_1$.
There is a first stage $s+1 > s_1$ where $i$ requires attention through (2); let $x_0 = source(s+1)$. 
If $y = cand(s+1) \not\in A$ then there is a stage $s' > s$ such that $s+1$ is the least $n > d_i$ with
$active(n) = 1$ and $E_{n-1} \subseteq A_s$. In the following stages when $i$ requires attention,
$y$ will be enumerated into $E_{s}, E_{s-1}, ...$, and finally into $E_{d_i}$, which gives the
desired contradiction. If $y \in A$, it might happen that $y$ is enumerated into $A$ before
it arrives in $E_{d_i}$. But then a new candidate $y'$ from $W_{g(x_0)}$ is chosen and a new attempt is
started to bring $y'$ into $E_{d_i}$. Again, it might happen that $y'$ is enumerated into $A$ before
it arrives in $E_{d_i}$. However, this process cannot repeat infinitely often, because otherwise
$x_0 \not\in B$ and hence there is some $y \in W_{g(x_0)}-A$. This $y$ would in some iteration
be chosen as a candidate which cannot be enumerated into $A$. So, at some stage $s_2+1 > s_1$
some $y$ is enumerated into $E_{d_i}$. Since $y \not\in A_{s_2}$ and $E_{d_i} \subseteq A_{s_2}$,
this implies that $E_{d_i}$ increases, a contradiction.

Thus, $i$ requires attention only finitely often and $(R_i)$ is satisfied. This completes the proof
of the inductive step.
\end{proof}

Recall that $A$ is strongly effectively simple 
if it is a coinfinite r.e. set and there is a total 
recursive function $f$ such that for all $e$,
$$W_e \subseteq \overline{A}\ \Rightarrow\ \max(W_e) < f(e).$$
Since every strongly effectively simple set is Q-complete 
\cite[Exercise III.6.21, a)]{Od89} we get the following corollary.

\begin{corollary} Every strongly effectively simple set has hard instances, even w.r.t.\ $\overline{ic}$.
\end{corollary}

It is known that hyperhypersimple sets are not Q-complete \cite[Theorem III.4.10]{Od89}, 
but we can still show that they have hard instances.

\begin{theorem} Every hyperhypersimple set has hard instances, even w.r.t.\ $\overline{ic}$.
\end{theorem}

\begin{proof} The basic idea of this proof is similar to the previous one.
Assume that $A$ is hyperhypersimple. We enumerate an r.e.\ sequence of finite 
sets $\{E_n\}_{n \in \nat}$ such that $|E_n| \leq 2^n$. It suffices to
get infinitely many $n$ such that there is $y \in E_n$ with $\overline{ic}(y:A) \geq n-2$.

Let $I_n = \{ p : l(p) < n-2\}$. We initialize $E_n = \{n\}$ and may later enumerate 
numbers from $E_n$ into $E_{n-1}$. This time we ensure that at any stage $s$ 
{\it at most two\/} numbers of $E_{n}$  belong to $\overline{A}_s$. 
We never enumerate a number twice into the same set.
Furthermore, we enumerate $x$ into $E_{n}$ at stage $s+1$
only if there is $p \in I_{n+1}$ which is $A$-consistent at stage $s+1$ and $U_s(p,x) = 0$.

From this invariant it already follows that $|E_n| \leq 2^n$: It is easy to see, by
induction on $k$, that we enumerate the $(2k+1)$-st number into $E_n$ at stage $s+1$
only if there are at least $k$ programs $p$ from $I_{n+1}$ which were $A$-consistent 
at some previous stage and are now diagonalized (i.e., for each such $p$ there is 
$z \in E_{n} \cap A_s$ such that $U_s(p,z) = 0$).  Since there are less than $2^{n-1}$
programs in $I_{n+1}$, it follows that $|E_n| < 2 \cdot 2^{n-1}+1 = 2^n+1$.

As in the previous proof, we say that $n$ is {\it saturated\/} at stage $s+1$ if
for every $y \in E_n$ there is $p \in I_n$ such that $p$ is $A$-consistent at 
stage $s+1$ and $U_s(p,y) = \chi_{A_s}(y)$. We want to produce infinitely many $n$
which are almost always not saturated.

To this end we construct for each $e$ a sequence $d^{e}_0 < d^{e}_1 <\cdots$ such that
for each $i$, $|\overline{A} \cap E_{d^e_i}| \geq 1$ or there is $n \in [d^e_i,d^e_{i+1})$
which is almost always not saturated. Suppose we have constructed at the end of stage $s$
an initial segment of this sequence, say $d^e_0<\cdots<d^e_{m+1}$. 
Let $count(n,s) = |\overline{A}_s \cap E_{n,s}|$.
We extend this initial segment at stage $s+1$ only if $count(d^e_i,s) \geq 1$
for all $i \leq m$. In the end we shall be able to argue that if the sequence is
infinite then there is a weak array which witnesses that $A$ is not hyperhypersimple.
Thus, the sequence must be finite, say $d^e_0<\cdots<d^e_{m(e)+1}$, and there is
$n \in [d^e_{m(e)}, d^e_{m(e)+1})$ which is almost always not saturated.
Also, since the strategy to extend the $e$-th sequence is active at only finitely
many stages, we can build an $(e+1)$-st sequence with $d^{e+1}_0 > d^e_{m(e)+1}$,
which will also be finite and gives us another number that is almost always not
saturated, etc.  

We assign priorities as follows: The definition of the $e$-th sequence has higher
priority than the definition of the $e'$-th sequence if $e<e'$. The definition
of the $i$-th member of the $e$-th sequence has higher priority than the definition
of the $i'$-th member if $i<i'$. Hence, we take the lexicographical ordering 
$<_{\mbox{\scriptsize lex}}$ on $\nat \times \nat$ as our priority ordering.

For technical reasons we enumerate for each $e$ a set $M_e$. When we are working on
the $e$-th sequence we try to establish for each $d^e_i$ a number $x \in E_{d^e_i} -A$.
In $M_e$ we enumerate the current candidate for $x$. 

We say that $(e,i)$ {\it requires attention\/} at stage $s+1$ if one of the following conditions
holds at the beginning of stage $s+1$.

\begin{enumerate}
\item[(1)] $d^e_i$ is undefined and for all $j \in [0,i-1)$:
$count(d^e_j,s) \geq 1$ and every $n \in [d^e_j, d^e_{j+1})$ is saturated at stage $s+1$.
\item[(2)] $d^e_i, d^e_{i+1}$ are both defined, $count(d^e_i,s) = 0$, and
every $n \in [d^e_i,d^e_{i+1})$ is saturated at stage $s+1$.
\end{enumerate}

\medskip
\noindent
{\it Construction:}\\
{\it Stage $0$}: Let $d^e_i = \uparrow$ and $M_e = \emptyset$ for all $e,i$, and  
let $E_n = \{n\}$ for all $n$.\\
{\it Stage $s+1$}: Choose the lexicographically least $(e,i)$ which requires attention
at stage $s+1$.

If it requires attention through (1) then let $d^e_i = s+1$, enumerate $s+1$ into $M_e$, 
and let $d^{e'}_j = \uparrow$ for all $(e',j) >_{\mbox{\scriptsize lex}} (e,i)$.

If it requires attention through (2) and there is a least $n \in (d^e_i,d^e_{i+1})$,
such that $count(n-1,s) \leq 1$ and there is a least $x \in E_{n,s}-(A_s \cup M_{e,s} \cup E_{n-1,s})$,
then enumerate $x$ into $E_{n-1}$. If in addition $n-1 = d^e_i$ then enumerate $x$ into $M_e$.
In any case, let $d^{e'}_j = \uparrow$ for all $(e',j) >_{\mbox{\scriptsize lex}} (e,i)$.\\
{\it End of Construction.}

\medskip
\noindent
It easily follows by induction on $s$ that $count(n,s) \leq 2$ for all $n,s$, in particular,
$|E \cap \overline{A}| \leq 2$. Also, we enumerate at stage $s+1$ a number $x$ from $E_n$ 
into $E_{n-1}$ only if it does not yet belong to $E_{n-1} \cap A$ and $n$ is saturated. In particular,
there is a program $p \in I_{n+1}$ which is $A$-consistent at stage $s+1$ and $U_s(p,x) = 0$.

\smallskip
\noindent
{\it Claim:\/} For every $e$, there are only finitely many stages where $(e,i)$
requires attention for some $i$.

\smallskip
\noindent
{\it Proof:\/}
Suppose for a contradiction that there exists a least $e$ and  infinitely many $s$ such that
$(e,i)$ requires attention at stage $s+1$ for some $i$. Then we argue that $A$ is not 
hyperhypersimple. First, there is a least stage $s_0 \geq 1$ such that no 
$(e',i')$ with $e' < e$ requires attention at any stage 
$s \geq s_0$. Then we define $d^e_0 = s_0$ at stage $s_0$ and we enumerate 
$s_0$ into $M_e$. By the choice of $s_0$, the value of $d^e_0$ has stabilized.
Note that all numbers which have been previously enumerated 
into $M_e$ are less than $s_0$ and so they do not matter for the following.
By induction on $s \geq s_0$, it follows that $E_{n,s}$ contains at most
one number from $M_{e,s}-A_s$ for all $n \geq d^e_0$.

Now we distinguish two cases:\\
(a) If there is a least $i$ such that $(e,i)$ requires attention infinitely often
then there is a stage $s_1 \geq s_0$ where all $d^e_j$ with $j \leq i$ have stabilized.
Thus, $(e,i)$ infinitely often requires attention through (2) and $d^e_{i+1}$ tends
to infinity. But then it follows similarly as in the previous proof that
unboundedly many numbers are eventually enumerated into $E_{d^e_i}$ which contradicts
the fact that the cardinality of $E_{d^e_i}$ is bounded: 

If $(e,i)$ requires attention through (2) at any stage $s \geq s_1$ then 
$count(d^e_i,s) = 0$, thus an $(e,j)$ with $j >i$ cannot require attention 
through (2) at any later stage $s' > s$, until a new number is enumerated into $E_{d^e_i}$ 
and $count(d^e_i,s') = 1$. During that time $M_e$ does not change. This guarantees that
eventually a new number is enumerated into $E_{d^e_i}$, since there exist numbers
$z \in (\bigcup_{n > d^e_i} E_{n,s}) - (A \cup M_e \cup E_{d^e_i,s})$. Since
$|E_{n,s} \cap (M_{e,s}-A_s)| \leq 1$ for $n \geq d^e_i$, it causes no problems
to maintain the constraint that a number $x$ is enumerated from $E_n$ into $E_{n-1}$ 
at stage $s+1$, only if $x\not\in (M_{e,s} \cup A_s)$.

\noindent
(b) If for every $i$ there are only finitely many  stages (but at least one stage)
where $(e,i)$ requires attention, then it follows that the values $d^e_i$ stabilize
and form an infinite increasing sequence. Let $d^e_i$ denote the final value.
Since the sequence is infinite it follows that $\lim_s count(d^e_i,s)  \geq 1$,
thus $|E_{d^e_i} \cap \overline{A}| \geq 1$. From the actual construction we get 
$|E_{d^e_i} \cap \overline{A}|=1$ and $E_{d^e_i} \cap \overline{A} \subseteq M_e$.

Uniformly in $i$ we enumerate an r.e.\ set $U_i$ as follows: If there is a stage
$s+1 \geq s_0$ where $(e,i)$ is the least pair which requires attention 
through (2) and a number $x$ is enumerated into $E_{d^e_i}$, then enumerate $x$
into $U_i$.

Since each such $x$ is also enumerated into $M_e$ and is therefore blocked
for the other sets,  it follows that the $U_i$'s are pairwise disjoint.
By the remarks above, each $U_i$ intersects $\overline{A}$. Thus, $A$ is not
hyperhypersimple. This contradiction completes the proof of the claim. $\Box$

\smallskip
\noindent
Thus, for each $e$ there exists a maximal $m(e) \geq 0$ such that the value
of $d^e_{m(e)+1}$ stabilizes and no $(e,j)$ with $j > m(e)+1$ requires
attention at any sufficiently large stage. This means that  there
exists $n \in [d^e_{m(e)},d^e_{m(e)+1})$ which is almost always not saturated.
Thus, there is $y \in E_n$ with $\overline{ic}(y:A) \geq n-2$. 
Clearly, we get infinitely many pairwise different such $y$'s.
This completes the proof.
\end{proof}

The previous result does not hold for hypersimple sets, since one can construct a hypersimple
set that does not have hard instances. This can be done, e.g., by a direct modification of the proof
of the next theorem.

Recall that $A$ is effectively simple if it is a coinfinite r.e.\ set and there is a
recursive function $f$ such that for all $e$,
$$W_e \subseteq \overline{A}\ \Rightarrow\ |W_e| \leq f(e).$$
It is known that every effectively simple set is T-complete \cite[Proposition III.2.18]{Od89}. 

\begin{theorem} There is an effectively simple set which does not have hard instances.
In particular, there is a T-complete set which does not have hard instances.
\end{theorem}

\begin{sketch} The construction in the proof of Theorem \ref{th:icc1} is not combinable
with the requirement of making $A$ effectively simple. Therefore, we use a modified
version were we do not attempt to have the instance complexity as low as possible.

In the following we outline the construction. $A$ will be effectively simple for some
$f$ to be determined later. As in the proof of Theorem \ref{th:icc1} we are given
a uniformly r.e.\ sequence $\{E_k\}_{k \in \nat}$ and we build  a partial
recursive function $\psi$ such that for almost all $k$ and for each 
$x \in E_k$ there is some $p \in \{0,1\}^k$ witnessing that $ic_{\psi}(x:A) \leq k$.

How do we define $\psi_p$? We will keep a list $S = S_k$ of programs of length $k$.
The length of $S$ will be fixed (depending on $k$). Furthermore, we have a pool $P = P_k$ of
unused programs of length $k$. At the beginning $|S| + |P| = 2^k$. During the construction
some of the programs in $S$ may become  inconsistent with $A$, in which case they are removed from
$S$ and new programs from $P$ are inserted into $S$. There may also exist a ``back-up 
program'' chosen from $P$.

The programs in $S$ will be defined at $x$ with a 0/1-value only if $x$ was enumerated 
into $E_k$. The definition proceeds in a round-robin fashion: 
The first program in $S$ takes care of the first number
which is enumerated into $E_k$, the second program takes care of the second number,
and so on. In this way we handle the first $|S|$ numbers. Ideally, we would
like that again the first program takes care of the $(|S|+1)$-st number, etc. However,
this does not work, because as soon as a program was brought into play
we have to define it for larger and larger inputs. So it might happen that
all of our programs are already defined (with output $\perp$) at $x$ when $x$
is enumerated as the $|S|+1$-st number at stage $s$.

Thus, we are using a program $q$ from $P$ which is still everywhere undefined and
define it as $\chi_{A_s}(z)$ for all $z<s$, in particular this covers 
all numbers currently in  $E_k$. For all larger values we
output $\perp$. $q$ is called the current back-up.
We also suspend defining the programs in $S$ until new numbers
$x \geq s$ are enumerated into $E_k$. Then we continue as above for the next
$|S|$ such numbers. After that a new program from $P$ is defined as the current
back-up in a similar way as $q$, and so on.

What is the advantage of that scheme? It is more robust against injuries
which may happen when a number $x$ with $\psi_p(x) = 0$ is later enumerated 
into $A$. In that case only one $p \in S$ is destroyed. Also, only the 
$x \in E_k$ are critical because for $x\not\in E_k$ we have $\psi_p(x) = \perp$.
If $p$ is destroyed then we assign a new program from $P$ as a substitute.

A crucial part in this process is the definition of the new back-up $q$
when a round has been completed at the beginning of stage $s$.
Before we define $\psi_q$, we enumerate all $x < s$ into $A$ which do not 
belong to any $E_n$ with $n < g(k)$: This defines the current $A_s$.
Here $g$ is some fast growing function to be determined later. 
Then we define $\psi_q(x) = \chi_{A_s}(x)$ for all $x<s$, 
and $\psi_q(x) = \perp$ otherwise.

We use the following strategy to make $A$
effectively simple. If at the end of some stage $s$ 
we have $W_{e,s} \subseteq \overline{A}_s$
and $|W_{e,s}| > f(e)$, then choose an $x \in W_{e,s}$ 
which does not belong to any $E_n$ with $n \leq g(e)$ 
and enumerate it into $A$. Note that $x$ exists if 
we choose $f$ large enough such that 
$f(e) \geq |E_0|+|E_1|+\cdots+|E_{g(e)}|$.

This completes the description of the construction.
It remains to choose the parameters such that it works.
We first count how many of the $\psi_p$ with $l(p) = k$ are used. 
Then we choose $|E_k|$ and $g$ in such a way 
that the number of used programs  is less than $2^k$.

Let $m = \max\{ n : g(n) \leq k \}$. Then for each $i \leq m$ there can be
$\lceil |E_i|/|S_i| \rceil$ many rounds and after each round all
programs in $S_k$ may be destroyed (and have to be replaced by new 
ones from $P_k$). At this time it is important
that after the action of $i$ we immediately define the new programs
that replace the former ones which have been destroyed.
We can do this without any further enumeration of elements into $A$.
There is no cascading effect which could blow up the number of 
injuries. Thus, at most $|S_k|\Sigma_{i=1}^m \lceil |E_i|/|S_i| \rceil$
many programs in $S_k$ are ever injured.

How many of the back-up functions are destroyed? Note that this may happen each time
when some $i<k$ acts, i.e., whenever $i$ completes a round. Thus, 
at most $\Sigma_{i=0}^{k-1} \lceil |E_i|/|S_i| \rceil$ many back-up functions
are destroyed. 

The number of injuries from making $A$ effectively simple
can be bounded by $m+k$: If we act for the sake of $W_{e,s} \cap A \not= \emptyset$
(which happens at most once), then a program from $S_k$ can be destroyed only if
$e < m$, and a current back-up program can be destroyed only if $e < k$. To see the
latter, note that if the current $q$ is defined at $x \not\in E_0\cup\ldots\cup E_{g(k)}$,
then $\psi_q(x) \in \{ 1, \perp\}$ because of the additional enumeration of numbers
into $A$ which was performed when $q$ was brought into play.

Thus, we need to ensure that for almost all $k$:
$$(+)\quad  2^k\ >\ |S_k| \Sigma_{i=1}^m \lceil |E_i|/|S_i| \rceil\ +\ 
\Sigma_{i=1}^{k-1} \lceil |E_i|/|S_i| \rceil\ + m + k.$$

Let $|S_k| = \lfloor 2^{k}/k \rfloor$, $g(k) = 2^k$, and
$E_k = \{ x : C(x) < 3k/2 \}$, so $|E_k| < 2^{3k/2}$.
Define the recursive function $f$ by $f(e) = \lceil \sum_{i=0}^{g(e)} 2^{3i/2} \rceil$.
The right hand side of $(+)$  is bounded above by
$$(2^{k}/k)(\log k)^2 \sqrt k\ +\ k^2 2^{k/2} + \log k + k$$
which is less than $2^k$ for all sufficiently large $k$.

With this choice of parameters we get for almost all $x$,
$C(x) \geq (3/2) (ic_{\psi}(x:A) - 1)$,  i.e., $ic_{\psi}(x:A) \leq (2/3) C(x)+1$.
Thus, $A$ does not have hard instances.
\end{sketch}

The previous results characterize the reducibilities $\leq_r$
with r~$\in  \{$m, btt, c, d, p, tt, wtt, Q, T$\}$ (cf.\ the figure in \cite[p.\ 341]{Od89}) 
such that every {r-complete} set has hard instances, for both $ic$ and $\overline{ic}$. 
In the following table we have marked the possible combinations.

\begin{center}
\renewcommand{\arraystretch}{1.5}
\begin{tabular}{||l||c|c|c|c|c|c|c|c|c|c|c|c||}         \hline
r  & m & btt & c & d & p & tt & wtt & Q & T \\  \hline\hline
$ic$ & $\times$ &  $\times$ &  $\times$ & $\times$ & $\times$ 
& $\times$ & $\times$ & $\times$ & \\ \hline
$\overline{ic}$ &  $\times$ & $\times$ & $\times$ & & & & & $\times$ & \\ \hline
\end{tabular}
\end{center}

\noindent
{\it Remark:\/} The T-degrees of r.e.\ sets with hard instances do not
coincide with any of the known degree classes. It can be shown that they 
form a proper subclass of the r.e.\ nonrecursive degrees and that they 
properly extend the array nonrecursive degrees.

\medskip
\noindent
{\bf Acknowledgments:} I would like to thank Lance Fortnow, Bill Gasarch,
and Paul Vit{\'a}nyi for comments on a preliminary version of this paper.

\end{document}